\documentclass[12pt]{amsart}

\usepackage[all]{xy}
\usepackage{xypic}
\usepackage{color}
\usepackage{graphicx}
\usepackage{amssymb, latexsym}

\usepackage[ german, frenchb, english]{babel}
\newtheorem{thm}{Theorem}[section]
\newtheorem{prop}[thm]{Proposition}
\newtheorem{cor}[thm]{Corollary}

\usepackage{float}
\usepackage{epsfig}
\restylefloat{table}
\usepackage{rotating}
\usepackage{chngcntr}
\counterwithout{table}{section}

\usepackage{colortbl}

\DeclareMathOperator{\elm}{elm}

\DeclareMathOperator{\Con}{{Con}}
\DeclareMathOperator{\Bun}{Bun}
\DeclareMathOperator{\Higgs}{Higgs}

\def\C{\mathbb{C}}
\def\P{\mathbb{P}}

\def\OO{\mathcal O}
\def\OOP{\OO_{\P^1}}

\def\p{\boldsymbol{p}}

\def\BUN{\mathfrak{Bun}}
\def\HIGGS{\mathfrak{Higgs}}

\usepackage{tikz}

\def\commutatif{\ar@{}[rd]|{{\LARGE \circlearrowright}}}

\numberwithin{equation}{section}

\begin{document}

\baselineskip=15.5pt

\title{Hitchin Hamiltonians in genus 2}

\author[V.Heu]{Viktoria Heu}
\address{IRMA, 7 rue Ren\'e-Descartes, 67084 Strasbourg Cedex, France}
\email{heu@math.unistra.fr}

\author[F.Loray]{Frank LORAY}
\address{IRMAR, Campus de Beaulieu, 35042 Rennes Cedex, France}
\email{frank.loray@univ-rennes1.fr}

\date{\today}


\thanks{\noindent The first author is supported by the ANR grants ANR-13-BS01-0001-01 and ANR-13-JS01-0002-01. \\ The second author is supported by CNRS}

\begin{abstract}
We give an explicit expression of the Hitchin Hamiltonian system for rank two vector bundles with trivial determinant bundle over a curve of genus two.
\end{abstract}

\maketitle

\section{Introduction}

We are interested in rank two vector bundles $E\to X$ with trivial determinant bundle $\det(E)=\mathcal{O}_X$ over a Riemann surface $X$ of genus $2$. The moduli space $\mathcal{M}_{NR}$ of semistable such vector bundles up to $S$-equivalence has been constructed by Narasimhan and Ramanan in \cite{NarasimhanRamanan}. If $E\in \mathcal{M}_{NR}$ is stable (and is therefore the unique vector bundle $S$-equivalent to $E$), the cotangent space of $\mathcal{M}_{NR}$ at $E$ is canonically isomorphic to the moduli space of trace free holomorphic Higgs fields $\theta: E\to E\otimes \Omega^1_X$ on $E$:
$$\mathrm{T}^\vee_E\mathcal{M}_{NR}\simeq \Higgs (X)|_E.$$
Since $  \mathcal{M}_{NR}\simeq \mathbb{P}^3$ and the locus of semistable but non stable bundles (up to $S$-equivalence) there is given by a singular quartic hypersurface, we have 
$$\mathrm{T}^\vee \mathcal{M}_{NR}\simeq \Higgs (X)$$ in restriction to a Zariski open subset of $\mathcal{M}_{NR}$, where $\Higgs (X)$ denotes the moduli space of tracefree holomorphic Higgs bundles $(E,\theta)$.

In \cite{Hitchin}, Hitchin considered the map \begin{equation}\label{HitchinMap} \mathrm{Hitch}:\left\{\begin{array}{ccc}  \Higgs (X)&\to &\mathrm{H}^0(X, \Omega^1_X\otimes \Omega^1_X)\\(E,\theta)&\mapsto & \det(\theta)\end{array}\right\}\end{equation} and established that it defines an algebraically completely integrable Hamiltonian system: the Liouville form on $\mathcal{M}_{NR}$ induces a symplectic structure on $\Higgs (X)$ and any set of (three) generators of  $\mathrm{H}^0(X, \Omega^1_X\otimes \Omega^1_X)$  commutes for the induced Poisson structure. Moreover, fibers of the Hitchin map are open sets of abelian varieties whose compactification is given by the Jacobian of the spectral curve. A broad field of applications has been deduced from the various algebraic and geometric properties of the Hitchin system and its generalizations since then.

Of course the Hitchin system in  \cite{Hitchin} is defined in a more general setting, but in the present paper, we focus on the special case as above (rank $2$ vector bundles with trivial determinant over curves of genus $2$) and announce results of a forthcoming paper \cite{FrankViktoria}:
 \begin{itemize}
 \item[$\bullet$] We describe the moduli space $\Bun(X)$ of (not necessarily semistable) vector bundles $E$ equivariant under the hyperelliptic involution $\iota$ on $X$. On the categorical quotient we have a birational morphism $$\Bun(X)\stackrel{\sim}{\dashrightarrow}\mathcal{M}_{NR}.$$
 \item[$\bullet$] It is well-known that there is no universal bundles on a Zariski-open subset of $\mathcal{M}_{NR}$. Yet from the dictionary between equivariant bundles and parabolic bundles on the quotient established in \cite{Biswas} (see also \cite{BiswasHeu}) we obtain a rational two-cover $$\Bun(X/\iota)\stackrel{2:1}{\dashrightarrow}\Bun(X)$$ and we construct a universal bundle over affine charts of $\Bun(X/\iota)$ which can   be identified with the universal bundle  in \cite{Michele} obtained from different methods. 
 \item[$\bullet$] We deduce a universal family of Higgs bundles on affine charts of $\Higgs(X/\iota):=\mathrm{T}^\vee \Bun(X/\iota)$. Note that in restriction to the stable locus, $\Higgs(X/\iota)\to \Bun(X/\iota)$ is a  principal $\mathbb{C}^3$-bundle.
\item[$\bullet$] These explicit universal families allow us to calculate the determinant map on $\Higgs(X/\iota)$ explicitly, which by construction factors through the Hitchin map. We deduce an explicit expression of the Hitchin map \eqref{HitchinMap} completing partial results in \cite{Emma}.
 \end{itemize}

\section{The Narasimhan-Ramanan moduli space $\mathcal{M}_{NR}$}

Let us first briefly recall the classical Narasimhan-Ramanan construction. Let $E\to X$ be a semistable rank 2 vector bundle with trivial determinant bundle over a Riemann surface $X$ of genus $2$. The subset $$C_E:=\{L\in \mathrm{Pic}^1(X) ~|~ \dim \mathrm{H}^0(X, E\otimes L) >0\}$$ of $\mathrm{Pic}^1(X)$ defines a divisor $D_E$ on $\mathrm{Pic}^1(X)$ which is linearly equivalent to the divisor $2\Theta$ on $\mathrm{Pic}^1(X)$, where $\Theta$ denotes the theta-divisor defined by the canonical embedding of $X$ in  $\mathrm{Pic}^1(X)$. In that way, we associate to the vector bundle $E$ an element $D_E$ of the Narasimhan-Ramanan moduli space $$\mathcal{M}_{NR}:=\P \mathrm{H}^0(\mathrm{Pic}^1(X), \mathcal{O}_{\mathrm{Pic}^1(X)}(2\Theta))\simeq \mathbb{P}^3.$$ For each smooth (analytic or algebraic)  family $\mathcal{E}\to X\times T$ of  
semistable rank 2 vector bundle with trivial determinant bundle on $X$, the Narasimhan-Ramanan classifying map 
$$T\to  \mathcal{M}_{NR}; \quad t\to D_{\mathcal{E}|_{X\times \{t\}}}$$ is an (analytic or algebraic) morphism. Moreover, if $D_E=D_{E'}$, then the vector bundles $E$ and $E'$ are $S$-equivalent:
\begin{itemize}
\item either $E$ is stable and then $E'=E$
\item or $E$ is strictly semistable (\emph{i.e.} semistable but not stable) and then there are line subbundles $L$ and $L'$ of degree $0$ of $E$ and $E$ respectively, such that $L=L'$ or $L=L'^{\otimes -1}$
\end{itemize}

The strictly semistable locus in $\mathcal{M}_{NR}$, which we shall denote by $\mathrm{Kum}(X)$ is thus defined by an embedding
$$\left\{\begin{array}{ccc}\mathrm{Pic}^0(X)/_{\iota}&\hookrightarrow& \mathcal{M}_{NR}\\L\ \mathrm{mod}\ \iota &\mapsto& D_{L\oplus L^{\otimes -1}}=L\cdot \Theta +\iota^*L\cdot \Theta\end{array}\right\}.$$
Note that if $\iota$ denotes the hyperelliptic involution on $X$, then $\iota^*L= L^{\otimes -1}$ by a classical argument that will be recalled in Section \ref{CoordsOnNR}. 

\subsection{Straightforward coordinates on $\mathcal{M}_{NR}$}\label{CoordsOnNR}
Any compact connected Riemann surface $X$ of genus $2$ can be embedded into $\P^1\times \P^1$ and is given, in a convenient affine chart, by an equation of the form
\begin{equation}\label{F}X\, : \, y^2=F(x)\quad \textrm{ with } \quad F(x)=x(x-1)(x-r)(x-s)(x-t).\end{equation}
The hyperelliptic involution on $X$ then writes $\iota : (x,y)\mapsto (x,-y)$ and the induced projection on the Riemann sphere $\P^1$ is given by $\pi: (x,y)\mapsto x.$ Denote by $$W:=\{w_0, w_1, w_r, w_s, w_t, w_\infty\}$$ the six Weierstrass points on $X$ invariant under the hyperelliptic involution, given by $w_i=(i,0)$ for $i\neq \infty$ and $w_\infty=(\infty, \infty)$. We will write $\underline W$ for $\pi(W)$.

Recall that the rational map
$$\left\{\begin{array}{ccc}X^{2}&\dashrightarrow & \mathrm{Pic}^2(X)\\\{P,Q\}&\mapsto& [P]+[Q]\end{array}\right\}$$ 
is surjective. More precisely, it is a blow-up of the canonical divisor $$K_X\sim [P]+[\iota (P)] \quad \textrm{ for all }†\quad P\in X.$$
Moreover, $\mathrm{Pic}^2(X) \simeq \mathrm{Pic}^1(X);\, D\mapsto D-[w_\infty]$ is an isomorphism. Global sections of $\mathcal{O}_{\mathrm{Pic}^1}(2\Theta)$ thus correspond bijectively to symmetric meromorphic functions on $X\times X$ with polar divisor at most $2\Delta +2\infty_1+2\infty_2$, where $\Delta:=\{(P,\iota(P)) ~|~P\in X\}$ and $\infty_i:=\{(P_1,P_2) ~|~P_i=w_\infty\}$. Any set of (four) generators of the vector space $\mathrm{H}^0(\mathrm{Pic}^1(X), \mathcal{O}_{\mathrm{Pic}^1}(2\Theta))$ can be expressed as such meromorphic functions, the simplest one being certainly the set $\{1, \mathrm{Sum}, \mathrm{Prod}, \mathrm{Diag}\}$ of functions in $(P_1,P_2)\in X\times X$ given, for  $P_i=(x_i,y_i)$, by
\begin{equation}\begin{array}{rll}1 : (P_1,P_2) & \mapsto & 1\vspace{.2cm}\\ \mathrm{Sum} : (P_1,P_2)& \mapsto &x_1+x_2\vspace{.2cm}\\ \mathrm{Prod} : (P_1,P_2)&\mapsto& x_1x_2,\\
 \mathrm{Diag} :(P_1,P_2)& \mapsto &\left(\frac{y_2-y_1}{x_2-x_1}\right)^2-(x_1+x_2)^3+(1+r+s+t)(x_1+x_2)^2 \hspace{.1cm}+\\&&+x_1x_2(x_1+x_2)-(r+s+t+rs+st+tr)(x_1+x_2)\end{array}\end{equation}
We obtain coordinates on $\mathcal{M}_{NR}=\P\mathrm{H}^0(\mathrm{Pic}^1(X), \mathcal{O}_{\mathrm{Pic}^1}(2\Theta))\simeq \P^3$, where we identify a point $(v_0:v_1:v_2:v_3)$ with the push-forward $D_E$ on $\mathrm{Pic}^1(X)$ of the zero-divisor of the meromorphic function 
\begin{equation}\label{expressDE} v_0\cdot 1+v_1\cdot\mathrm{Sum}+v_2\cdot \mathrm{Prod}+v_3 \mathrm{Diag}\end{equation} on $X^{(2)}$. 

For example if $E=L\oplus L^{-1}$, where $$L=\mathcal{O}_X([Q_1]+[Q_2]-2[w_\infty])=\mathcal{O}_X([Q_1]-[\iota(Q_2)])$$
and $Q_1,Q_2\in X$, we can calculate the (unique up to scalar) meromorphic function as in \eqref{expressDE} whose $0$-divisor corresponds to $D_E=L\cdot \Theta +\iota^*L\cdot \Theta$ and we obtain
\begin{equation}\label{KumParam}(v_0:v_1:v_2:v_3)=(-\mathrm{Diag}(Q_1,Q_2) : \mathrm{Prod}(Q_1,Q_2) : -\mathrm{Sum}(Q_1,Q_2) : 1).\end{equation}

The strictly semistable locus $\mathrm{Kum}(X)$ in $\mathcal{M}_{NR}$ is parametrized by $\left\{Q_1,Q_2\right\}\in X^{(2)}$ according to formula \eqref{KumParam}. We deduce an equation for $\mathrm{Kum}(X)$ in our coordinates $(v_0:v_1:v_2:v_3)$ of $\mathcal{M}_{NR}$ : 
\begin{equation}\label{KummerV}\begin{array}{rrcc}\mathrm{Kum}\left(X\right)~:\\
0=&(v_0v_2-v_1^2)^2&\cdot &1\vspace{.2cm}\\
&-2\left[[(\sigma_1+\sigma_2)v_1+(\sigma_2+\sigma_3)v_2](v_0v_2-v_1^2)\right.\\&+\left.2(v_0+\sigma_1 v_1)(v_0
+v_1)v_1+2(\sigma_2v_1+\sigma_3 v_2)(v_1+v_2)v_1\right]&\cdot&v_3\vspace{.2cm}\\
&-2\sigma_3(v_0v_2-v_1^2)+\left[\left[(\sigma_1+\sigma_2)^2v_1+(\sigma_2+\sigma_3)^2v_2\right](v_1+v_2)\right.\\&\left.-(\sigma_1+\sigma_3)^2v_1v_2+4[(\sigma_2+\sigma_3)v_0-\sigma_3v_2]v_1\right]&\cdot&\vspace{.2cm}v_3^2\\
&-2\sigma_3\left[(\sigma_1+\sigma_2)v_1-(\sigma_2+\sigma_3)v_2\right]&\cdot&\vspace{.2cm}v_3^3\\
&+ \sigma_3^2&\cdot&v_3^4.\end{array}\end{equation}
Since the  strictly semistable locus $\mathrm{Kum}(X)$ is a quartic with sixteen conic singularities it is usually referred to as \emph{the Kummer surface} in the context of $\mathcal{M}_{NR}$.

Let $ \mathcal{O}_X(\tau)$ be a $2$-torsion line bundle on $X$, \emph{i.e.} $\mathcal{O}_X(2\tau)\simeq \mathcal{O}_X$. Then
\begin{equation}\label{tau}\tau\sim [w_i]-[w_j]\quad \textrm{ with }†w_i,w_j\in W \end{equation}
and the group of $2$-torsion line bundles on $X$ with respect to the tensor product is isomorphic to $\left(\mathbb{Z}/2\mathbb{Z}\right)^4$. 
 If $E$ is a rank two vector bundle with trivial determinant bundle over $X$, then its \emph{twist} $E\otimes \mathcal{O}_X(\tau)$ also has the trivial determinant line bundle.
Moreover, by construction of  the Narasimhan-Ramanan moduli space, the action of the group of $2$-torsion line bundles by twist is linear and free on $\mathcal{M}_{NR}$ and preserves $\mathrm{Kum}(X).$  By formula \eqref{KumParam} we can explicitly calculate the coordinates $(v_0:v_1:v_2:v_3)$ of the trivial bundle $E_0=\mathcal{O}_X\oplus \mathcal{O}_X$ and its twists $$E_\tau:=E_0\otimes \mathcal{O}_X(\tau).$$ The trivial bundle for example is given by $E_0 : (1:0:0:0)$.  
Note that these sixteen bundles $E_\tau$ correspond to the sixteen singularities of the Kummer surface $\mathrm{Kum}(X)$. 

The fact that we know the action (by permutation) of the $2$-torsion group on the set of bundles $E_\tau$ and we also know the coordinates of these bundles in the Narasimhan-Ramanan moduli space is sufficient to calculate explicitly the linear action of the $2$-torsion group on $\mathcal{M}_{NR}$ : for any $\tau$ as in \eqref{tau}, there is a matrix $M_\tau \in \mathrm{SL}_4\C$ such that if the image  of $E$ under the Narasimhan-Ramanan classifying map is given by $(v_0:v_1:v_2:v_3)$, then $E'=E\otimes \mathcal{O}_X(\tau)$ is given by $(v_0':v_1':v_2':v_3')$ with
$$\begin{pmatrix}v_0'\\v_1'\\v_2'\\v_3'
\end{pmatrix}=M_\tau \cdot \begin{pmatrix}v_0\\v_1\\v_2\\v_3
\end{pmatrix}.$$
The equivalence classes in $\mathrm{PGL}_{4}\C$ of these matrices (with respect to a set of generators of the $2$-torsion group) then are given by the following:
$$\begin{array}{ccl}M_{[w_0]-[w_\infty]}&\sim&\begin{pmatrix}0&rs+st+rt+rst& rst& 0\\ 0& 0& 0&rst\\ 1& 0& 0& -(rs+st+rt+rst)\\ 0& 1& 0&0
\end{pmatrix}\\
\\
M_{[w_1]-[w_\infty]}&\sim&\begin{pmatrix}1& r+s+t+rst& rs+st+rt& 0\\ -1& -1& 0& rs+st+rt\\ 1& 0& -1& -(r+s+t+rst)\\ 0& 1& 1& 1\end{pmatrix}
\\
\\M_{[w_r]-[w_\infty]}&\sim&\begin{pmatrix}r^2& r^2(1+s+t)+st& r^2(s+t+st)& 0\\-r& -r^2& 0& r^2(s+t+st)\\ 1& 0& -r^2& -r^2(1+s+t)-st\\ 0& 1& r& r^2
\end{pmatrix}\\
\\
M_{[w_s]-[w_\infty]}&\sim&\begin{pmatrix}s^2& s^2(1+r+t)+rt& s^2(r+t+rt)& 0\\ -s& -s^2& 0& s^2(r+t+rt)\\ 1& 0& -s^2& -s^2(1+r+t)-rt\\ 0& 1& s& s^2\end{pmatrix}
\end{array}$$

\subsection{Nice coordinates on $\mathcal{M}_{NR}$}
A quick calculation shows that the character of the representation 
$$\rho : \left\{\begin{array}{ccc} \left(\mathbb{Z}/2\mathbb{Z}\right)^4 & \to & \mathrm{SL}_4\C\\
\mathcal{O}_X(\tau) & \mapsto &  M_\tau \end{array}\right\}$$ 
introduced above is the regular one : it vanishes on all elements of the group exept $\mathcal{O}_X$.
Hence $\rho$ is conjugated for example to the regular representation $\widetilde{\rho} : \mathcal{O}_X(\tau)\mapsto \widetilde{M}_\tau$ given by
$$\begin{array}{ccccccccc}\widetilde{M}_{[w_0]-[w_\infty]}&=&\begin{pmatrix}0&0&1&0\\0&0&0&1\\1&0&0&0\\0&1&0&0
\end{pmatrix}&&, &&
\widetilde{M}_{[w_1]-[w_\infty]}&=&\begin{pmatrix}0&1&0&0\\-1&0&0&0\\0&0&0&-1\\0&0&1&0\end{pmatrix}
\\
\\
\widetilde{M}_{[w_r]-[w_\infty]}&=&\begin{pmatrix}1&0&0&0\\0&-1&0&0\\0&0&-1&0\\0&0&0&1
\end{pmatrix}&&, &&
\widetilde{M}_{[w_s]-[w_\infty]}&=&\begin{pmatrix}0&1&0&0\\1&0&0&0\\0&0&0&-1\\0&0&-1&0\end{pmatrix}
\end{array}
$$
Any conjugation matrix $M\in \mathrm{SL}_4\C$ such that $\widetilde{\rho}=M\rho M^{-1}$ is given, up to a scalar, as follows:
Choose square-roots $\omega_0, \omega_1, \omega_r, \omega_s$ such that
 $$
\omega_0^2=F'(0),\quad  \omega_1^2=-F'(1), \quad
\omega_r^2=F'(r),\quad   \omega_s^2=F'(s)
,$$
where $F(x)$ is given in \eqref{F} and $F'(x)$ is its derivative with respect to $x$.
 Then  
\begin{equation}\label{nicecoords}M={\begin{pmatrix}a&b&c&d\\-b&a&d&-c\\ c&d&a&b\\d&-c&-b&a\end{pmatrix}}\cdot\begin{pmatrix}1&1&0&-\omega_0\\0& \omega_1&0&0\\0&\omega_0&\omega_0&\omega_0\\0&0&0&\omega_0 \omega_1\end{pmatrix}, \vspace{.2cm}\end{equation}
 $\begin{array}{lccrcl}
\textrm{where }&&&a&=&rst(r-s) \omega_1+t\omega_r\omega_s\vspace{.2cm}-rt(r-1)\omega_s-st \omega_1\omega_r\vspace{.2cm}\\
&&&b&=&-st(s-1)\omega_r+rt \omega_1\omega_s\vspace{.2cm}\\
&&&c&=&t(r-s)\omega_0 \omega_1-t(r-1)\omega_0\omega_s\vspace{.2cm}\\
&&&d&=&-t(r-1)(s-1)(r-s)\omega_0+t(s-1)\omega_0\omega_r .\end{array}$\\
After the coordinate-change $(v_0:v_1:v_2:v_3)\mapsto (u_0:u_1:u_2:u_3)$ on $\mathcal{M}_{NR}$ defined by 
$${\begin{pmatrix}u_0\\u_1\\u_2\\u_3\end{pmatrix}}=M\cdot{\begin{pmatrix}v_0\\v_1\\v_2\\v_3\end{pmatrix}},$$
the action of the $2$-torsion group is then normalized to $\widetilde{\rho}$. In particular, the equation of the Kummer surface with respect to the coordinates $(u_0:u_1:u_2:u_3)$ is invariant under double-transpositions and double-changes of signs. Calculation shows
\begin{equation}\label{KummerU}\begin{array}{rrcc}\mathrm{Kum}\left(X\right)~:\\
0=& (u_0^4+u_1^4+u_2^4+u_3^4)-8\frac{rs-rt+r-s}{t(s-1)}u_0u_1u_2u_3-2\frac{st+t-2s}{t(s-1)}(u_0^2u_3^2+u_1^2u_2^2)\\&-2\frac{2r-t}{t}(u_1^2u_3^2+u_0^2u_2^2)+2\frac{2r-s-1}{s-1}(u_2^2u_3^2+u_0^2u_1^2).\end{array}\end{equation}

In summary, the straightforward coordinates $(v_0:v_1:v_2:v_3)$ of $\mathcal{M}_{NR}$ introduced in the previous section have the advantage that 
\begin{itemize}
\item[$\bullet$] a given divisor $D_E$ on $\mathrm{Pic}^1(X)$ linearly equivalent to $2\Theta$ can rather easily be expressed in terms of $(v_0:v_1:v_2:v_3)$,
\item[$\bullet$] and we are going to use this property when we describe the universal family on a $2$-cover of $\mathrm{M}_{NR}$,
\end{itemize}
whereas the new coordinates $(u_0:u_1:u_2:u_3)$ of $\mathcal{M}_{NR}$ defined above have the advantages that 
\begin{itemize}
\item[$\bullet$] the action of the $2$-torsion group is simply expressed by double-transpositions and double-changes of signs of $(u_0:u_1:u_2:u_3)$ and
\item[$\bullet$] the equation of the Kummer surface is rather symmetric. As pointed out in \cite{Emma}, the classical line geometry for Kummer surfaces in $\P^3$ is related to certain symmetries of the Hitchin Hamiltonians. For this geometrical reason, the explicit Hitchin Hamiltonians we are going to establish have a much simpler expression with respect to (dual) coordinates $(u_0:u_1:u_2:u_3)$ when compared to $(v_0:v_1:v_2:v_3).
$\item[$\bullet$] Moreover, the five $u_i$-polynomials in \eqref{KummerU} invariant under the action of the $2$-torsion group define a natural map $\mathcal{M}_{NR}\to \P^4$. The image is a quartic hyper surface \cite[Proposition 10.2.7]{Dolgachev} and can be seen as the coarse moduli space of semistable $\P^1$-bundles over $X$.
\end{itemize}

\section{How to construct a bundle from a point in $\mathcal{M}_{NR}$}
Given a stable rank 2 vector bundle $E$ with trivial determinant bundle on $X$, the Narasimhan-Ramanan divisor $D_E\in |2\Theta|$ can be seen as space of line subbundles $L$ of $E$ of degree $-1$. Whilst we know that for any $D\in |2\Theta|$ there is a semistable vector bundle $E$ with $D_E=D$, it is not obvious how to construct it. We provide such a construction by considering the moduli space of rank 2 vector bundles $E$ with trivial determinant bundle on $X$ equivariant under the hyperelliptic involution. For the present exposition however, we restrict our attention to the space
$\Bun(X)$
of rank 2 vector bundles $E$ with trivial determinant bundle on $X$ such that
\begin{itemize} 
\item[$\bullet$] $E$ is stable but \emph{off the odd Gunning planes}, which means that no line subbundle $L\subset E$ is isomorphic to $\mathcal{O}_X(-[w_i])$ for some $w_i\in W$, or
\item[$\bullet$] $E$ is strictly semistable but undecomposable, or
\item[$\bullet$] $E=L\oplus \iota^*L$ where $L=\mathcal{O}_X([P]-[Q])$ satisfies $P,Q \not \in W$, or
\item[$\bullet$] $E$ is an odd Gunning bundle, \emph{i.e} given by the unique non-trivial extension 
$$0\longrightarrow \mathcal{O}_X([w_i])\longrightarrow E \longrightarrow \mathcal{O}_X(-[w_i])\longrightarrow 0$$
for a Weierstrass point $w_i\in W$.
\end{itemize}
We construct $\Bun(X)$ as an algebraic stack whose categorical quotient is birational to the Narasimhan-Ramanan moduli space
$$ \Bun(X) \stackrel{1:1}{\dashrightarrow} \mathcal{M}_{NR}.$$ 

For convenience of notation let us for now denote by $ \Bun(X)$ the set of vector bundles $E$ as in the above list, before we put an algebraic structure on $\Bun(X)$.  We use the fact that any bundle $E\in \Bun(X)$ is equivariant under the hyperelliptic involution: 
\begin{prop}\label{h} Let $E$ be a vector bundle in $ \Bun(X)$. Then there is a bundle isomorphism $h$ such that the following diagram commutes
$$\begin{xy}\xymatrix{
E \ar[r]^\sim_h\ar@/_2pc/[rr]_{\mathrm{id}_E}&\iota^*E\ar[r]^\sim_{\iota^*h}&E .
}
\end{xy}$$
and such that for each Weierstrass point $w_i\in W$, the induced automorphism of the Weierstrass fibre
$$h|_{E_{w_i}}: E_{w_i}\to \iota^*E_{w_i}\simeq E_{w_i}$$
possesses two opposite eigenvalues $+1$ and $-1$. 
\end{prop}

Now, hyperelliptic descent \cite{Biswas}
$$\pi_*E=(\underline E^+, \underline \p^+)\oplus (\underline E^-, \underline \p^-)$$
produces two rank $2$ vector bundles $\underline E^\pm$ with determinant bundle $\det(\underline E^\pm)=\mathcal{O}_{\P^1}(-3)$ over the Riemann sphere, each endowed with a natural quasi-parabolic structure $\underline \p^\pm$ with support $\underline W=\pi(W)$. Moreover,

\begin{prop}\label{understandLift} 
Consider $(E,h)$ as in Proposition \ref{h}. Denote by $\p^+$ and $\p^-$ the quasi-parabolic structure with support $W$ on $E$ induced by the $+1$ and $-1$ eigendirections of $h$ respectively. Then
$$(E, \p^{\pm})= \elm_W^+( \pi^*(\underline E^\pm, \underline \p^\pm)),$$ where $ \elm_W^+$ denotes the composition of six positive elementary transformations, one over each Weierstrass point, given by the corresponding quasi-parabolic direction of $ \pi^*(\underline E^\pm, \underline \p^\pm)$.
\end{prop}

In convenient local coordinates $(\zeta ,Y)\in U\times \C^2$ of $E\to X$ near a Weierstrass point $w_i:\{\zeta=0\}$, the map $\elm_W^+\circ \pi^*$ can be understood as follows: 

$$\begin{xy}
\xymatrix{
  \p|_{w_i} :  \{Y\in \mathrm{Vect}_\C {\left(\begin{smallmatrix}1\\0\end{smallmatrix}\right) }\}&\left( \zeta, Y\right) \ar[r]^h& \left( -\zeta, {\left(\begin{smallmatrix}1&0\\0&-1\end{smallmatrix}\right) }Y\right) \\
  \hat{\p}|_{w_i} :  \{\widehat{Y}\in  \mathrm{Vect}_\C {\left(\begin{smallmatrix}0\\1\end{smallmatrix}\right) } \}& \left( \zeta, \widehat{ Y}\right) =\left(\zeta,{\left(\begin{smallmatrix}1&0\\0&\frac{1}{\zeta}\end{smallmatrix}\right) }Y\right)\ar[r]^h\ar[u]^{\elm^+_{\hat{\p}}}& \left( -\zeta, \widehat{ Y}\right)\ar[u]^{\elm^+_{\hat{\p}}}\\
  \underline{\p}|_{\pi(w_i)} : \{ \underline{Y}\in \mathrm{Vect}_\C {\left(\begin{smallmatrix}0\\1\end{smallmatrix}\right) } \}& \left( \underline \zeta, \underline{ Y}\right) =\left(\zeta^2,\widehat{Y}\right)\ar[r]^{\mathrm{id}}\ar[u]^{\pi^*}& \left( \underline\zeta, \underline{ Y}\right)\ar[u]^{\pi^*}\\
}
\end{xy}$$

Let $\mu$ be a real number in $[0,1]$.
Denote by $\Bun_{\mu}(X/\iota)$ the moduli space of pairs $(\underline E , \underline \p)$, where $\underline E$ is a rank 2 vector bundle of degree $-3$ over $\P^1$ and $\underline \p$ is a quasi-parabolic structure with support $\underline W$ such that
$(\underline E, \underline \p)$ is a stable parabolic bundle if to each quasi-parabolic direction $\underline \p|_{\underline w_i}$ we associate the parabolic weight $\mu$. 
For each choice of $\mu$, this moduli space is either empty or birational to $\P^3$ \cite{LoraySaito}. Moreover, for any $\mu\in [0,1]$, the map  $$\mathcal{O}_{\P^1}(-3)\otimes \mathrm{elm}^+_{\underline W}$$ is a canonical birational isomorphism between $\Bun_{\mu}(X/\iota)$ and 
$\Bun_{1-\mu}(X/\iota).$

Note further that for $\mu=\frac{1}{5}$, the space $\Bun_{\mu}(X/\iota)$ is precisely the moduli space of those quasi-parabolic bundles $(\underline E, \underline \p)$, where $\underline E$ is a vector bundle on $\P^1$ and $\underline \p$ is a quasi parabolic structure with support $\underline{W}$ on $\underline E$ such that
\begin{itemize}
\item[$\bullet$]$\underline E = \mathcal{O}_{\P^1}(-1)\oplus \mathcal{O}_{\P^1}(-2)$,
\item[$\bullet$] the quasi-parabolic directions $\underline{\p}_{\underline w_i}$ are all disjoint from the total space of the destabilizing subbundle $\mathcal{O}_{\P^1}(-1)\subset \underline E$ and
\item[$\bullet$] the quasi-parabolic directions $\underline{\p}_{\underline w_i}$ are not all contained in the total space of a same subbundle $\mathcal{O}_{\P^1}(-2)\hookrightarrow \underline E$
\end{itemize}

Consider the following affine chart $(R,S,T)\in \mathbb{C}^3$ of $\Bun_{\frac{1}{5}}(X/\iota)$, which we shall call \emph{the canonical chart} : Recall that $\underline E = \mathcal{O}_{\P^1}(-1)\oplus \mathcal{O}_{\P^1}(-2)$.
Let $\sigma_1$ be a meromorphic section of some line subbundle $\mathcal{O}_{\P^1}(-2)\hookrightarrow \underline E$ with only one (double) pole over $x=\infty$. Let $\sigma_{-1}$ be a meromorphic section of the unique line subbundle $\mathcal{O}_{\P^1}(-1)\subset \underline E$ with only one pole over $x=\infty$. In the total space of $\underline E$ restricted to $\P^1\setminus \{\infty\}$, we consider coordinates  $\left( x, \left(\begin{smallmatrix} z_1\\z_2\end{smallmatrix}\right)\right)$ given by $(x,z_1\sigma_{-1}+z_2\sigma_{1})$. 
To $(R,S,T)\in \mathbb{C}^3$ we then associate the following normalized quasi-parabolic structure on $\underline E$:
\begin{equation}\label{RiemScheme}\begin{matrix}x=0 & x=1 & x= r & x=s & x=t & x=\infty\\
\begin{pmatrix}0\\ 1\end{pmatrix} &\begin{pmatrix}1\\ 1\end{pmatrix} &\begin{pmatrix}R\\ 1\end{pmatrix} &\begin{pmatrix}S\\ 1\end{pmatrix} &\begin{pmatrix}T\\ 1\end{pmatrix} & \OOP(-1)
\end{matrix}.\end{equation}

Here the first line indicates  the Weierstrass point $\underline w_i$ we are considering, whereas the second line defines a generator of the corresponding quasi-parabolic direction. Note that \eqref{RiemScheme} already defines a universal quasi-parabolic bundle over the canonical chart.  The lifting map $\elm_W^+\circ \pi^*$ is well-defined and algebraic and provides a universal rank $2$ vector bundle with trivial determinant bundle over the canonical chart of $\Bun_{\frac{1}{5}}(X/\iota)$.

\begin{prop}\label{PropFormulesRSTtoNR}
The Narasimhan-Ramanan  classifying map $\C^3_{(R,S,T)}\dashrightarrow\mathcal{M}_{NR}$ is explicitely 
given by $(R,S,T)\mapsto(v_0:v_1:v_2:v_3)$ where
 $$\begin{array}{rcl}
 v_0&=& s^2t^2(r^2-1)(s-t)R-r^2t^2(s^2-1)(r-t)S+s^2r^2(t^2-1)(r-s)T+\\ &&+t^2(t-1)(r^2-s^2)RS-s^2(s-1)(r^2-t^2)RT+r^2(r-1)(s^2-t^2)ST\vspace{.3cm}\\
 v_1& =&rst\left[( (r-1)(s-t)R-(s-1)(r-t)S+(t-1)(r-s)T+\right.\\ &&\left.+ (t-1)(r-s)RS-(s-1)(r-t)RT+(r-1)(s-t)ST \right]\vspace{.3cm}\\
 v_2&=&-st(r^2-1)(s-t)R+rt(s^2-1)(r-t)S-rs(t^2-1)(r-s)T-\\&&-t(t-1)(r^2-s^2)RS+s(s-1)(r^2-t^2)RT-r(r-1)(s^2-t^2)ST\vspace{.3cm}\\
 v_3&=& st(r-1)(s-t)R-rt(s-1)(r-t)S+sr(t-1)(r-s)T+\\
 &&+ t(t-1)(r-s)RS-s(s-1)(r-t)RT+r(r-1)(s-t)ST
 \end{array}$$
 The indeterminacy points
 $$(R,S,T)=(0,0,0),\ \ \ (1,1,1) \ \ \ \text{and}\ \ \ (r,s,t)$$ of this map 
 correspond to the odd Gunning bundles $E_{[w_1]}, E_{[w_0]}$ and $E_{[w_\infty]}$ respectively.
Conversely, a generic point  $(v_0:v_1:v_2:v_3)\in \mathcal{M}_{NR}$ has precisely two preimages in $\C^3_{(R,S,T)}$ given by
$$\begin{array}{rcl}
R&=&\frac{r(t-1)(v_0+rv_1-r(s+t+st)v_3)T}{t(r-1)(v_0+tv_1-t(r+s+rs)v_3)-(r-t)(v_0+v_1-(rs+st+rt)v_3)T}\vspace{.3cm}\\
S&=&\frac{s(t-1)(v_0+sv_1-s(r+t+rt)v_3)T}{t(s-1)(v_0+tv_1-t(r+s+rs)v_3)-(s-t)(v_0+v_1-(rs+st+rt)v_3)T} ,
\end{array}$$
where $T$ is any solution of $aT^2+btT+ct^2=0$ with $$\begin{array}{rcl}
a&=&(v_1+v_2t+v_3t^2)(v_0+v_1-(rs+st+rt)v_3)\\
b&=&-(1+t)(v_0v_2+v_1^2+tv_1v_3)-2(v_0v_1+tv_0v_3+tv_1v_2)\\&&+(rs+st+rt)(tv_1+v_2+tv_3)v_3+(r+s+rs)(v_1+t^2v_2+t^2v_3)v_3\\
c&=&(v_1+v_2+v_3)(v_0+tv_1-t(r+s+rs)v_3). \end{array}$$ The discriminant of this polynomial leads again to the equation \eqref{KummerV} of the Kummer surface.
\end{prop}

 By construction, $\Bun_{\frac{1}{5}}(X/\iota)$ is covered by affine charts similar to the canonical chart, where we just permute the role of the Weierstrass points in \eqref{RiemScheme}. The (birational) transition maps between affine charts are obvious. The Galois-involution $\mathcal{O}_{\P^1}(-3)\otimes \mathrm{elm}^+_{\underline W}$ 
 is given in the canonical chart  by the birational map $(R,S,T) \mapsto(\widetilde{R},\widetilde{S},\widetilde{T})$, where $$\begin{array}{l}
 \widetilde{R}= \lambda(R,S,T)\cdot  \frac{(s-t)+(t-1)S-(s-1)T}{-t(s-1)S+s(t-1)T+(s-t)ST}\\
\widetilde{S}= \lambda(R,S,T)\cdot \frac{(r-t)+(t-1)R-(r-1)T}{-t(r-1)R+r(t-1)T+(r-t)RT}\vspace{.3cm}\\
\widetilde{T}=\lambda(R,S,T)\cdot \frac{(r-s)+(s-1)R-(r-1)S}{-s(r-1)R+r(s-1)S+(r-s)RS}\vspace{.3cm}
\end{array}$$
$$\text{and}\ \ \ \begin{array}{l}\lambda(R,S,T)=\frac{t(r-s)RS-s(r-t)RT+r(s-t)ST}{(s-t)R-(r-t)S+(r-s)T}\end{array}.$$  
For the Galois-involution to be everywhere well defined, we need to consider the smooth (non separated) scheme $\Bun(X/\iota)$ obtained by canonically gluing 
$\Bun_{\frac{1}{5}}(X/\iota)$ and $\Bun_{\frac{4}{5}}(X/\iota)$.
From an exhaustive case-by case study,  one can show that $\Bun(X)$ corresponds precisely to the isomorphism classes of the Galois-involution on $\Bun(X/\iota)$. In terms of parabolic bundles, the Galois involution is given by $\elm_W^+\circ \pi^*$. In other words,

\begin{prop} The map $\elm_W^+\circ \pi^* : \Bun(X/\iota)\stackrel{2:1}{\longrightarrow} \Bun(X)$ is an algebraic $2$ cover. \end{prop}

Moreover, the lift of the Kummer surface in $\mathcal{M}_{NR}$ defines a dual Weddle surface in $\Bun_{\frac{1}{5}}(X/\iota)\simeq \P^3 \subset \Bun(X/\iota)$ which is given with respect to the canonical chart by the equation 
$$\begin{array}{rrcl}\mathrm{Wed}(X) &:&\\
0&=& ((s-t)R+(t-r)S+(r-s)T)RST+t((r-1)S-(s-1)R)RS\\
&&+r((s-1)T-(t-1)S)ST+s((t-1)R-(r-1)T)RT\\
&&-t(r-s)RS-r(s-t)ST-s(t-r)RT.
\end{array}$$

\section{Application to Higgs bundles}

Let $E$ again be a rank 2 vector bundle over $X$. By definition, the moduli space of tracefree Higgs fields on $E$ is given by $\mathrm{H}^0(X, \mathfrak{sl}(E)\otimes  \Omega^1_X)$, where $\mathfrak{sl}(E)$ denotes the vector bundle of trace-free endomorphisms of $E$. By Serre duality, we have 
$$\mathrm{H}^0(X, \mathfrak{sl}(E)\otimes  \Omega^1_X)\simeq \mathrm{H}^1(X, \mathfrak{sl}(E)^\vee \otimes  \Omega^1_X)^\vee.$$
If $\det(E)=\mathcal{O}_X$, then $\mathfrak{sl}(E)^\vee = \mathfrak{sl}(E)$. If $E$ is stable, then $\mathfrak{sl}(E)$
possesses no non-trivial global sections and then $\mathrm{H}^1(X, \mathfrak{sl}(E)^\vee \otimes  \Omega^1_X)$ canonically identifies with 
$\Higgs (X):=T^\vee\mathcal{M}_{NR}$. We can calculate explicitly the Hitchin map $$ \mathrm{Hitch}:\left\{\begin{array}{ccc}  \Higgs (X)&\to &\mathrm{H}^0(X, \Omega^1_X\otimes \Omega^1_X)\\(E,\theta)&\mapsto & \det(\theta)\end{array}\right\}$$ from the following idea :  The complement in $\Bun(X)$ of the image of the Weddle surface is embedded into $\mathcal{M}\setminus \mathrm{Kum}(X)$ (we obtain all stable bundles except those on the odd Gunning planes).
Since we have a universal vector bundle in each affine chart of the two-cover $\Bun(X/\iota)$ of $ \Bun(X)$, we can expect to find a universal family of Higgs bundles there as well.  Then we calculate a Hitchin map for $\Bun(X/\iota)$ and push it down to $\mathcal{M}_{NR}$. 

More precisely, we will calculate the Hitchin map in the following steps :
\begin{itemize}
\item[$\bullet$] Provided that $\mathrm{H}^0(\P^1,\mathfrak{sl} (\underline{E}, \underline{\p}))=\{0\}$, we have a canonical isomorphism  $$\mathrm{T}_{(\underline E , \underline p)}\Bun (X/\iota) = \mathrm{H}^1(\P^1,\mathfrak{sl} (\underline{E}, \underline{\p})),$$ 
  where $\mathfrak{sl} (\underline E, \underline \p))$ denotes trace free endomorphisms of $\underline E$ leaving $\underline \p$ invariant. 
 We  work out how the hyperelliptic descent $\phi : \mathrm{elm}^+\circ \pi^*$ defines an algebraic $2$-cover
$$\Higgs (X/\iota):= \mathrm{T}^\vee\Bun (X/\iota)  \stackrel{\phi}{\longrightarrow} \mathrm{T}^\vee\Bun (X) $$
\item[$\bullet$] The Liouville form on $\Bun(X/\iota)$ is given with respect to coordinates $(R,S,T)$ of the canonical chart by 
$\mathrm{d}R+\mathrm{d}S+\mathrm{d}T.$
We work out Serre duality for the generators $$\frac{\partial}{\partial R} \,, \frac{\partial}{\partial S}\,, \frac{\partial}{\partial T} \in \mathrm{T} \Bun(X/\iota)$$ and deduce an explicit universal Higgs bundle on an affine chart of $\Higgs (X/\iota)$. 
\item[$\bullet$] We calculate the determinant map $$\Higgs (X/\iota)\longrightarrow \mathrm{H}^0(\P^1, \Omega^1_{\P^1}\otimes \Omega^1_{\P^1}(\underline W))\simeq \mathrm{H}^0(X, \Omega^1_{X}\otimes \Omega^1_{X})$$ and show that it factors through the Hitchin map. 
Then we deduce the explicit Hitchin map from the formulas in Proposition \ref{PropFormulesRSTtoNR}.\end{itemize}

\subsection{Hyperelliptic descent, again}
 One can show that if $E$ is a stable rank 2 vector bundle with trivial determinant on $X$ and $h$ is a lift of the hyperlliptic involution as in Proposition \ref{h}, then any trace free Higgs field $\theta$ on $E$ is $h$-equivariant, that is, the following diagram commutes :
$$\begin{xy}\xymatrix{
E \ar[r]^\theta \ar[d]^h& E\otimes \Omega^1_X\\
\iota^*E\ar[r]^{\iota^*\theta}&\iota^*E\otimes \Omega^1_X\ar[u]^{\iota^*h} .
}
\end{xy}$$
Hyperelliptic decent of the pair $(E,\theta)$ then produces two triples $(\underline E , \underline \theta , \underline p)$, 
where $$\underline \theta : \underline E \to \underline E \otimes \Omega^1_{\P^1}(\underline W)$$ 
is a logarithmic Higgs field with at most apparent singularities over $\underline W$ : the residue at any $\underline w_i \in \underline W$ is either zero or conjugated to $$\mathrm{Res}_{\underline{w}_i}(\underline \theta) \sim \begin{pmatrix} 0&0\\1&0\end{pmatrix}$$
such that the quasi-parabolic $\underline{\p}_{\underline{w}_i}$ corresponds precisely to the $0$-eigendirection of the residue.
Conversely, consider a logarithmic Higgs field $$\underline{ \widetilde{\theta}} \in \mathrm{H}^0(\P^1 , \mathfrak{sl}(\underline E, \underline {\widetilde{\p}} ) \otimes \Omega^1_{\P^1}(\underline W))$$ that  lies in the image of the canonical embedding 
$$\mathrm{H}^0(\P^1,\mathfrak{sl} (\underline E, \underline \p)\otimes\Omega^1_{\P^1}) \hookrightarrow \mathrm{H}^0(\P^1,\mathfrak{sl} (\underline E, \underline \p)\otimes\Omega^1_{\P^1}(\underline W )).$$ In other words, $\underline{  \widetilde{\theta} }$ has only trivial residues:  $$\mathrm{Res}_{\underline{w}_i}(\widetilde{\underline{\theta}}) \sim \begin{pmatrix} 0&0\\0&0\end{pmatrix}.$$
Then the logarithmic Higgs field 
 $\underline \theta$ obtained from  $\widetilde{\underline{\theta}}$ by applying the meromorphic gauge transformation $$\OO_{\P^1}(-3)\otimes\mathrm{elm}_{\underline{W}}^+$$ 
 has (at most) apparent singularities over each Weierstrass point, and the quasi-parabolic structure $\underline{\p}$ obtained from $\underline{\widetilde{\p}}$ corresponds to  $0$-eigendirections of $\underline \theta$. Note that here we have to choose a meromorphic section $\sigma : \P^1\to \OO_{\P^1}(-3)$.
As long as we are not on the Kummer, respectively Weddle surface, we obtain an isomorphism
$$\begin{xy}\xymatrix{\mathrm{T}^\vee_{(\underline E ,\widetilde{\underline{\p}})}\Bun (X/\iota)=\mathrm{H}^1(\P^1,\mathfrak{sl} (\underline E , \widetilde{\underline{\p}}))^\vee\ar[r]_\sim^{\mathrm{Serre}}\ar[dd]_\sim&\mathrm{H}^0(\P^1,\mathfrak{sl} (\underline E , \widetilde{\underline{\p}})\otimes\Omega^1_{\P^1})\ar[d]_\sim^{\mathcal{O}_{\P^1}(-3)\otimes\elm^+_{\underline{W}}}\\ 
& \mathrm{H}^0(\P^1,\mathfrak{sl} (\underline E\,, \underline{\p})\otimes\Omega^1_{\P^1}(\underline W))^{\mathrm{apparent } }\ar[d]_{\sim}^{ \elm^+_{W}\circ \pi^*}\\\mathrm{T}^\vee_{E}\Bun (X)=\mathrm{H}^1(X,\mathfrak{sl} (E))^\vee&\mathrm{H}^0(X,\mathfrak{sl} (E)\otimes\Omega^1_{X})\ar[l]_\sim^{\mathrm{Serre}}
}\end{xy}$$
and deduce an algebraic $2$-cover $$\mathrm{Higgs}(X/\iota)\stackrel{2:1}{\longrightarrow}\mathrm{Higgs}(X).$$

\subsection{Universal Higgs bundles} Serre duality  gives us a perfect pairing
$$\langle\cdot,\cdot\rangle: \left\{\begin{array}{ccc} \mathrm{H}^1(\P^1,\mathfrak{sl}(\underline{E},\underline{\p}))\times \mathrm{H}^0(\P^1,\mathfrak{sl}_2(\underline E\, , \underline \p)\otimes\Omega^1_{\P^1}(\underline W ))^{\mathrm{apparent}}&\to&\C\\
(\phi  , \underline \theta) & \mapsto  &\sum \mathrm{Res}( \mathrm{trace}(\phi\cdot  \underline \theta))
\end{array}\right\}.$$
Let $(\underline E , \underline \p)$ be an element of $\Bun (X/\iota)$ given with respect to the canonical chart by $(R_0,S_0,T_0)\in \mathbb{C}^3$. 
The vector field $\frac{\partial}{\partial R}\in \mathrm{T}_{(R_0,S_0,T_0)}\Bun (X/\iota)$ is given in $\mathrm{H}^1(\P^1 , \mathfrak{sl}(\underline E , \underline \p))$  by the cocycle
$$\phi_{01}:=\begin{pmatrix} 0 & 1 \\0& 0\end{pmatrix}$$
with respect to trivialization charts 
$U_0\times \C^2$ with $U_0:=\P^1 \setminus \{r\}$ and  $U_1\times \C^2$ with $U_1:=D_{\varepsilon}(r)$ of $\mathfrak{sl}(\underline E)$. 
Indeed, if we consider $\exp(\zeta \phi)=\left(\begin{smallmatrix} 1 & \zeta \\0& 1\end{smallmatrix}\right)$ as applied from the left over $U_0\cap U_1\subset U_1$ to $U_0\cap U_1\subset U_0$ for a quasi-parabolic vector bundle $(\underline E , \underline \p)$ with parabolic structure normalized as in \eqref{RiemScheme}, we obtain the quasi parabolic structure corresponding to $(R+\zeta, S, T)$. 
The dual basis in $ \mathrm{H}^0(\P^1,\mathfrak{sl}_2(\underline E\, , \underline \p)\otimes\Omega^1_{\P^1}(\underline W ))^{\mathrm{apparent}}$ with respect to $\langle\cdot,\cdot\rangle$ of the basis 
$$\left(\frac{\partial}{\partial R}, \frac{\partial}{\partial S}, \frac{\partial}{\partial T}\right)$$
of $\mathrm{T}_{(R_0,S_0,T_0)}\Bun (X/\iota)$ then is given by $\left( \underline{\theta}_r, \underline{\theta}_s, \underline{\theta}_t\right)$ with 
$$\begin{array}{rcl}\underline{\theta}_r&:=&\begin{pmatrix}0&0\\ 1-R&0\end{pmatrix}\frac{dx}{x}
+\begin{pmatrix}R& -R\\ R&-R\end{pmatrix}\frac{dx}{x-1}
+\begin{pmatrix} -R & R^2\\ -1& R\end{pmatrix}\frac{dx}{x-r}\vspace{.3cm}\\\underline{\theta}_s&:=&\begin{pmatrix}0&0\\ 1-S&0\end{pmatrix}\frac{dx}{x}
+\begin{pmatrix}S& -S\\ S&-S\end{pmatrix}\frac{dx}{x-1}
+\begin{pmatrix} -S & S^2\\ -1& S\end{pmatrix}\frac{dx}{x-s}\vspace{.3cm}\\
\underline{\theta}_t&:=&\begin{pmatrix}0&0\\ 1-T&0\end{pmatrix}\frac{dx}{x}
+\begin{pmatrix}T& -T\\ T&-T\end{pmatrix}\frac{dx}{x-1}
+\begin{pmatrix} -T & T^2\\ -1& T\end{pmatrix}\frac{dx}{x-t} .\end{array}$$

We obtain the universal Higgs bundle $\widetilde{\underline \theta}=\mathcal{O}_{\P^1}(-3)\otimes \mathrm{elm}^+_{\underline W}(\theta)$ defined by \begin{equation}\label{universalHiggs}\underline{\theta}=c_r\underline{\theta}_r+c_s\underline{\theta}_s+c_t\underline{\theta}_t\end{equation} on the canonical chart $(R,S,T,c_r,c_s,c_t)\in \mathbb{C}^6$ of $\mathrm{T}^\vee\Bun(X/\iota)$. Recall that all other charts are obtained up to Galois involution by permuting the role of the Weierstrass points.

\begin{cor}The Liouville form on $\Bun(X/\iota)$, given with respect to the canonical chart $(R,S,T)\in \mathbb{C}^3$ by
$$\mathrm{d}R+\mathrm{d}S+\mathrm{d}T$$ defines a holomorphic 
symplectic $2$-form on $\Higgs(X / \iota)$ given with respect to the canonical chart $(R,S,T,c_r,c_s,c_t)\in \mathbb{C}^6$ of $\mathrm{T}^\vee\Bun(X/\iota)$ by
$$\mathrm{d}R\wedge \mathrm{d}c_r+\mathrm{d}S\wedge \mathrm{d}c_s+\mathrm{d}T\wedge \mathrm{d}c_t.$$
\end{cor}

\subsection{The Hitchin fibration}
The determinant map 
\begin{equation}\label{detmap}\left\{\begin{array}{ccc}\mathbb{C}^6& \to& \mathrm{H}^0(\P^1, \Omega^1_{\P^1}\otimes \Omega^1_{\P^1}(-6))\\
(R,S,T,c_r,c_s,c_t) & \mapsto & \det (\underline \theta)=(h_2x^2+h_1x+h_0)\frac{(dx)^{\otimes2}}{x(x-1)(x-r)(x-s)(x-t)}\end{array}\right\},\end{equation}
where $\underline{\theta}$ is the Higgs bundle in \eqref{universalHiggs} is given by 

\begin{table}[H]
$$\begin{array}{rcl}
h_0&=& \left( c_r(R-1)+c_s(S-1)+c_t(T-1)
 \right)  \left(c_rst (R-r)R+c_srt(S-s)S+c_trs(T-t)T\right)\\\\
h_1&=&+{c_r}\,
 \left(c_r(s+t)(r+1)+c_ss(t+1)+c_tt(s+1) \right) {R}^{2}-{{\it c_r}}^{2} \left( {\it t}+{\it s} \right) {R}^
{3}\\&&+{\it c_s}
\, \left(c_s(r+t)(s+1)+c_rr(t+1)+c_tt(r+1)
 \right) {S}^{2}-{{\it c_s}}^{2} \left( {\it t}+{\it r} \right) {S}^
{3}\\&&+{\it c_t}\, \left( {\it c_t}(r+s)(t+1)+{\it c_r}r(s+1)+{\it c_s}{\it s}(r+1)
 \right) {T}^{2}-{{\it c_t}}^{2} \left( {\it r}+{\it s} \right) {T}^{
3}\\&&-c_rc_s(t(R-1+S-1)+r(S-s)+s(R-r))RS \\&&-c_rc_t(s(R-1+T-1)+r(T-t)+t(R-r))RT\\&&-c_sc_t(r(S-1+T-1)+s(T-t)+t(S-s))ST\\&&- \left( {\it c_t}{\it t}(r+s)+{\it c_r}{\it r}({\it s}+t) +{
\it c_s}{\it s}(r+t) \right) (c_rR+c_sS+c_tT)
\\
\\
h_2&=&\left(c_r(R-1)R+c_s(S-1)S+c_t(T-1)T\right)  \left( {\it c_r}(R-r)+{\it c_s}(S-s)+
{\it c_t}(T-t) \right)
\end{array}
$$
\caption{Explicit Hitchin Hamiltonians for the canonical coordinates $(R,S,T)$ on  $\Bun(X/\iota)$}\label{HitchinRST}
\end{table}

It is easy to check that the Hitchin Hamiltonians $h_0,h_1,h_2$ do Poisson-commute as expected :
for any $f,g\in\{h_0,h_1,h_2\}$, we have
$$\sum_{i=r,s,t}\frac{\partial f}{\partial p_i}\frac{\partial g}{\partial q_i}-\frac{\partial f}{\partial q_i}\frac{\partial g}{\partial p_i}=0$$
in Darboux notation $(p_r,p_s,p_t,q_r,q_s,q_t):=(R,S,T,c_r,c_s,c_t)$.

Since the determinant is invariant under meromorphic gauge transformations (and in particular elementary transformations), we can immediately deduce the Hitchin map 
$$\underline{\mathrm{Hitch}}: \left\{\begin{array}{ccc}\mathrm{T}^\vee\Bun(X/\iota)& \to& \mathrm{H}^0(\P^1, \Omega^1_{\P^1}\otimes \Omega^1_{\P^1})\\
(R,S,T,c_r,c_s,c_t) & \mapsto & \det (\underline{\widetilde \theta})= (h_2x^2+h_1x+h_0)\frac{(dx)^{\otimes2}}{\sigma^2\cdot x(x-1)(x-r)(x-s)(x-t)}\end{array}\right\},$$
where  $\sigma$ is our previously chosen meromorphic section $\sigma: \P^1\to \OO_{\P^1}(-3)$ and $\underline{\widetilde \theta}$ is the universal Higgs bundle in the canonical chart of $\mathrm{T}^\vee\Bun(X/\iota)$. 

More importantly, again since the determinant map does is invariant under meromorphic Gauge transformations, the map in \eqref{detmap} 
factors through the Hitchin map $\HIGGS(X)\simeq T^\vee\BUN(X) \to \mathrm{H}^0(X, \Omega^1_X\otimes \Omega^1_X)$ by construction. 
Consider the natural rational map $\phi^*:\mathrm{T}^\vee \mathcal{M}_{NR}\dashrightarrow \mathrm{T}^\vee \Bun(X/\iota)$
induced by the map $\phi:\Bun(X/\iota)\dashrightarrow\mathcal{M}_{NR}$ stated explicitly with respect to the canonical chart in Proposition \ref{PropFormulesRSTtoNR}. The general section $c_r\mathrm{d}R+c_s\mathrm{d}S+c_t\mathrm{d}T$ then lifts to a general section $\mu_0 \mathrm{d}\left(\frac{v_0}{v_3}\right)+\mu_1 \mathrm{d}\left(\frac{v_1}{v_3}\right)+\mu_2 \mathrm{d}\left(\frac{v_2}{v_3}\right)$. Moreover, from the explicit coordinate change $(v_0:v_1:v_2:v_3)\leftrightarrow (u_0:u_1:u_2:u_3)$ to the nice coordinates, given in \eqref{nicecoords}, we know how to identify general sections
  $$\eta_0 \mathrm{d}\left(\frac{u_0}{u_3}\right)+\eta_1 \mathrm{d}\left(\frac{u_1}{u_3}\right)+\eta_2 \mathrm{d}\left(\frac{u_2}{u_3}\right)= \mu_0 \mathrm{d}\left(\frac{v_0}{v_3}\right)+\mu_1 \mathrm{d}\left(\frac{v_1}{v_3}\right)+\mu_2 \mathrm{d}\left(\frac{v_2}{v_3}\right).$$
  The Hamiltonians $h_0,h_1,h_2$ of the Hitchin map on $\mathcal{M}_{NR}$ then can be explicitly deduced from \eqref{detmap}. We get 
$${\mathrm{Hitch}}: \left\{\begin{array}{ccc}\mathrm{T}^\vee\mathcal{M}_{NR}& \to& \mathrm{H}^0(X, \Omega^1_{X}\otimes \Omega^1_{X})\\
((u_0:u_1:u_2:u_3), \eta_0,\eta_1,\eta_2) & \mapsto &  (h_2x^2+h_1x+h_0)\frac{(dx)^{\otimes2}}{x(x-1)(x-r)(x-s)(x-t)}\end{array}\right\},$$
where
\vspace{-.5cm}\\
\bgroup
\def\arraystretch{1.3}
\begin{table}[H]
$$\begin{array}{l}
\begin{array}{rcl}
h_0&=&  \frac{1}{4 u_3^4} \cdot \left\{ 
\begin{array}{rl}
rst\cdot & \left[\eta_0 (u_0^2-u_3^2)+\eta_1 (u_0 u_1+u_2 u_3)+\eta_2 (u_0 u_2+u_1 u_3)\right]^2\\
-st\cdot & \left[\eta_0 (u_0 u_1-u_2 u_3)+\eta_1 (u_1^2+u_3^2)+\eta_2 (u_0 u_3+u_1 u_2)\right]^2\\
+4rs\cdot & \left(\eta_0 u_0+\eta_1 u_1\right)^2u_3^2\\
-rt \cdot & \left[\eta_0 (u_0^2+u_3^2)+\eta_1 (u_0 u_1+u_2 u_3)+\eta_2 (u_0 u_2-u_1 u_3)\right]^2
\end{array}\right.
\end{array}
\vspace{.2cm}\\

\begin{array}{rcl}
h_1&=&  \frac{1}{4 u_3^4} \cdot \left\{ 
\begin{array}{rl}
t\cdot &  \left(u_0^2+u_1^2+u_2^2+u_3^2\right)  \left[(\eta_0^2+\eta_1^2+\eta_2^2) u_3^2+(\eta_0 u_0+\eta_1 u_1+\eta_2 u_2)^2\right]\\
+st\cdot &  \left(u_0^2-u_1^2+u_2^2-u_3^2\right)  \left[(\eta_0^2-\eta_1^2+\eta_2^2) u_3^2-(\eta_0 u_0+\eta_1 u_1+\eta_2 u_2)^2\right]
\\
+4r\cdot &  \left(u_0 u_2-u_1 u_3\right) u_3  \left[\eta_0 \eta_2 u_3+(\eta_0 u_0+\eta_1 u_1+\eta_2 u_2) \eta_1\right]
\\
+4sr\cdot &  \left(u_0 u_2+u_1 u_3\right) u_3  \left[\eta_0 \eta_2 u_3-(\eta_0 u_0+\eta_1 u_1+\eta_2 u_2) \eta_1\right]
\\
+4s\cdot &  \left(u_0 u_3+u_1 u_2\right) u_3  \left[\eta_1 \eta_2 u_3-(\eta_0 u_0+\eta_1 u_1+\eta_2 u_2) \eta_0\right]
\\
+4rt\cdot &  \left(u_0 u_1+u_2 u_3\right) u_3  \left[\eta_0 \eta_1 u_3-(\eta_0 u_0+\eta_1 u_1+\eta_2 u_2) \eta_2\right]
\end{array}\right.
\end{array}

\vspace{.2cm}\\

\begin{array}{rcl}
h_2&=&  \frac{1}{4 u_3^4} \cdot \left\{ 
\begin{array}{rl}
s\cdot &  \left[\eta_0 (u_0 u_2+u_1 u_3)+\eta_1 (u_0 u_3+u_1 u_2)+\eta_2 (u_2^2-u_3^2)\right]^2
\\
-1\cdot & \left[\eta_0 (u_0 u_2-u_1 u_3)+\eta_1 (u_0 u_3+u_1 u_2)+\eta_2 (u_2^2+u_3^2)\right]^2
\\
-t\cdot &  \left[\eta_0 (u_0 u_1+u_3 u_3)-\eta_2 (u_0 u_3-u_1 u_2)+\eta_1 (u_2^2+u_3^2)\right]^2\\
+4 r \cdot &  \left(\eta_1 u_1+\eta_2 u_2\right)^2 u_3^2

\end{array}\right.
\end{array}
\end{array}$$
\caption{Explicit Hitchin Hamiltonians for the coordinates $(u_0:u_1:u_2:u_3)$ of $\mathcal{M}_{NR}$. 
}\label{HitchinT}
\end{table}

Note that in \cite{Emma}, B. van Geemen and E. Previato conjectured a projective version of explicit Hitchin Hamiltonians, which has been confirmed in \cite{GTNB}. These Hamiltonians $H_1,\ldots H_6$ can be seen as evaluations, up to functions in the base, of the explicit Hitchin map at the Weierstrass points. More precisely, if we denote $$h(x):=h_2x^2+h_1x+h_0,$$ where $h_i$ for $i\in \{0,1,2\}$ then
$$\begin{array}{rcc cc rcc}
H_1 &=& \frac{4h(0)}{rst}  &\textrm{ }  &\textrm{ } & H_4  &=& \frac{4h(s)}{s(s-1)(s-r)(s-t)}\vspace{.3cm}\\
H_2 &=& -\frac{4h(t)}{t(t-1)(t-r)(t-s)} &&& H_5 &=& \frac{4h(r)}{r(r-1)(r-s)(r-t)}\vspace{.3cm}\\
H_3 &=& \frac{4h(1)}{(r-1)(s-1)(t-1)}&&& H_6 &=&0 .
\end{array}$$


\end{document}